\documentclass[letterpaper, 10 pt, conference]{ieeeconf}  

\IEEEoverridecommandlockouts                              
\usepackage[all]{xy}
\usepackage{color}
\usepackage{subfigure}
\usepackage{bm}
\usepackage{graphicx}
\usepackage{color}
\usepackage{wrapfig}
\usepackage{epsf}
\usepackage{times,mathptm}
\usepackage{url}
\usepackage{amssymb, amsmath}

\title{\LARGE\bf
DistFlow ODE: Modeling, Analyzing and Controlling\\
Long Distribution Feeder
}

\author{Danhua Wang, Konstantin Turitsyn, and Michael Chertkov
\thanks{
The work of MC at LANL was carried out under the auspices of the
National Nuclear Security Administration of the U.S. Department of
Energy at Los Alamos National Laboratory under Contract No.
DE-AC52-06NA25396. The work of MC was funded in part by DTRA/DOD under the grant BRCALL06-Per3-D-2-0022 on ``Network Adaptability from WMD Disruption and Cascading Failures".
The work of KT and MC is also partially supported by the collaborative (and linked) NSF grants at MIT and NMC on ``Spectroscopy of Power Grids".}
\thanks{D. Wang is with Department of Mathematics of SMU, Dallas, TX  75275,
{\tt\small danhuaw@mail.smu.edu}
}
\thanks{K. Turitsyn is with Department of Mechanical Engineering at MIT,
77 Massachusetts Ave., Cambridge MA 02139, {\tt\small turitsyn@mit.edu}}
\thanks{M. Chertkov is with Theory Division \& Center for Nonlinear Studies at LANL,
Los Alamos, NM 87545 and also with New Mexico Consortium, Los Alamos, NM 87544 {\tt\small chertkov@lanl.gov}}
}

\begin{document}

\maketitle
\thispagestyle{empty}
\pagestyle{empty}

\begin{abstract}
We consider a linear feeder connecting multiple distributed loads and generators to the sub-station. Voltage is controlled directly at the sub-station, however, voltage down the line shifts up or down, in particular depending on if the feeder operates in the power export regime or power import regime. Starting from this finite element description of the feeder, assuming that the consumption/generation is distributed heterogeneously along the feeder, and following the asymptotic homogenization approach, we derive simple low-parametric ODE model of the feeder. We also explain how the homogeneous ODE modeling is generalized to account for other distributed effects, e.g. for inverter based and voltage dependent control of reactive power. The resulting system of the DistFlow-ODEs, relating homogenized voltage to flows of real and reactive power along the lines, admits computationally efficient analysis in terms of the minimal number of the feeder line ``media" parameters,  such as the ratio of the inductance-to-resistance densities. Exploring the space of the media and control parameters allows us to test and juxtapose different measures of the system performance, in particular expressed in terms of the voltage drop along the feeder, power import/export from the feeder line as the whole, power losses within the feeder, and critical (with respect to possible voltage collapse) length of the feeder.
Our most surprising funding relates to performance  of a feeder rich on PhotoVoltaic (PV) systems  during a sunny day.  We observe that if the feeder is sufficiently long the DistFlow-ODEs may have multiple stable solutions. The multiplicity may mean troubles for successful recovery of the feeder after a very short,  few periods long, fault at the head of the line. 
\end{abstract}

\section{Introduction}
\label{sec:Intro}
{\bf Motivation.}
Typical feeder of a distribution system forms a line bringing power to a large number of sequentially-connected customers. Electric characteristics of the feeder are normally uniform while a snapshot of consumption varies from neighbor to neighbor even if these are statistically similar customers. We are interested to describe power flows and voltage profile along the line, where the latter is maintained constant only at the head of the line. When feeder consists of standard consumers, which extract both real and reactive power from the system, the voltage sages down the line showing the largest deviation from the nominal value at the end of the line.  Standard design principe for this type of long but simple line is to have it limited in length,  so that the typical voltage drop at the end of the line would not exceed the 5\% which can still be tolerated. However,  the consumption pattern vary  thus leading from time to time to an incidental off-bounds drop or raise of voltage. One may also discuss lines of varying length,  asking, in particular, what is the maximal possible length of a feeder which supports (within given model of consumption) a valid power flow solution? Moreover,  one is interested to generalize, thus discussing not only feeders consuming electricity but also producing, i.e. these built of distributed generators, e.g. representing small-scale PV-systems.  Such PV-enforced customers are consumers when it is cloudy turning into producers on a sunny day \cite{lopes2007integrating}. Line consisting of distributed producers may inject power in the transmission system (if real losses along the line are small), also leading to voltage increase at the end of the line with respect to the level set at the head of the line.

{\bf Related Work.} Discussion of the static Power Flow (PF) equations in the case of a distribution feeder was pioneered in \cite{89BWa,89BWb}. It was shown,  that it is natural and convenient to state the PF equations in this one-dimensional case in terms of voltages, and real and reactive powers flowing through segments of the line.  In this study we will also adopt the so-called DistFlow formulation of the PF equations from \cite{89BWa,89BWb}. The main technical thread of the present manuscript is related to
the DistFlow equations considered in the limit of large (continuously many) consumers forming the feeder. In particular we derive and analyze the Distflow ODE version of the original DistFlow. ODE of this type was first introduced and analyzed in \cite{11CBTCL},  where the main focus was on recovering the classic power engineering phenomena of voltage collapse \cite{Weedy1968,Venikov1975,Venikov1977,Taylor1994,Cutsem1998,VanCutsem2000} in the continuous setting of a long statistically homogeneous feeder. Note also that a related description of the power system as a media was made in \cite{Parashar2004,Thorp1998},  where propagation of electro-mechanical waves through transmission system was introduced and discussed. (This type of waves do not develop in the transmission system, simply because of the lack of big inertial machines in the distribution.  Also, and for the same reason, the electro-mechanical waves are not generated within the feeder.) In addition to following the Distflow ODE thread of \cite{11CBTCL},  this paper also continues  and extends analysis and approach of \cite{10TSBCa,10TSBCb,11TSBC} in what concerns inverter-based distributed voltage control in the distribution feeder.

{\bf Our Contribution.} The task of this manuscript is two fold.  First, and taking the algorithmic/computational stance, we aim at developing low-parametric (spatially homogenized) ODE modeling of a feeder consisting of many similar consumers. In some of the most general regimes, e.g. in the regime correspondent to constant consumption/production of power along the feeder, ODE model posses symmetry which allows funding power flows and voltage profile most efficiently. On the other hand,  and once the homogenized ODE model is stated and justified,  we aim to analyze power flows in the long distributed feeder,  paying a special attention to discovering qualitative phenomena of interest for operation and control of the feeder in the interesting regimes. To complete the second task we analyzed distributed feeder in different regimes.  For example,  we study the feeder in the regime of a standard distributed consumption of real power and analyze different voltage control schemes associated with the freedom in choosing consumption or injection of reactive power through the distributed inverters.  We also analyze future possibility of using feeder rich on PV systems in the
regime where the feeder would provide export of power to the transmission (high voltage) part of the grid.  Our most surprising conclusions apply to this case.  We have found that this export regime may be extremely dangerous as generating more than one stable solutions.  Then the healthy (normal voltage) coexists with possibly many low-voltage branches. Normally the healthy branch is to be realized,  however a significant perturbation,  e.g. a short fault at the head of the line, may force the system to transition into a low-voltage branch where the system would stay till the inverters are disconnected/tripped (or burned). The dangerous states are complex,  in particular showing oscillations of power flow directions along the feeder.  We also find that some previously suggested reactive control,  even though improving the quality of the healthy solution,  does not remove the dangerous degeneracy.

The material in the rest of the manuscript is organized as follows. We state the standard DistFlow formulation and show how the homogenization approach results in the DistFlow ODEs in Section \ref{sec:Hom}. We review recently suggested reactive flow/voltage control schemes in Section \ref{sec:Control}. Section \ref{sec:Re-scaling} is devoted to description of how the ODE formulation of the power flow equations allows re-scaling and even further reduction in the number of key characteristics and consequently efficient computations of the Cauchi (initial value) type. We present simulations of the feeder working in the consumption/generation,  control/no-control regimes and discuss the results in Section \ref{sec:Simulations}. Section \ref{sec:Conclusions} summarizes the material and suggests a path forward.

\section{Homogenization: from Finite-Element DistFlow to DistFlow-ODEs}
\label{sec:Hom}
The DistFlow equations of Baran, Wu \cite{89BWa,89BWb} have the following form
\begin{eqnarray}
&& P_{k+1}-P_k=p_k-r_k \frac{P_k^2+Q_k^2}{v_k^2},
\label{DF1}\\
&& Q_{k+1}-Q_k=q_k-x_k\frac{P_k^2+v_k^2}{v_k^2},
\label{DF2}\\
&& v_{k+1}^2-v_{k}^2=-2 (r_kP_k+x_k Q_k)-(r^2_k+x^2_k)\frac{P_k^2+Q_k^2}{v_k^2},
\label{DF3}
\end{eqnarray}
where $k=0,\cdots, N-1$ enumerates buses of the feeder, sequentially connected in a line, and
$P_k,Q_k$ stand for real and reactive power flowing from bus $k$ to bus $k+1$. $v_k$ is voltage, while $p_k$ and $q_k$ describe the overall consumption (negative/positive) of real and reactive powers at the bus $k$. The values of $r_k$ and $x_k$ represent the resistance and inductance of the line element connecting $k$ and $k+1$ buses. Finite feeder line with the tap changing transformer in the beginning can be modeled by the boundary conditions in the rescaled p.u. variables:
\begin{eqnarray}
v_0 = 1,\quad P_{N+1} = Q_{N+1} = 0.
\label{DF-BC}
\end{eqnarray}
These boundary conditions close the system of equations (\ref{DF1},\ref{DF2},\ref{DF3}) and allow full reconstruction of the power and voltage profile given the consumption/injection, $p_k, q_k$.

When the line is long and the number of consumers is large,  $N\gg 1$,  one is looking to simplify the system of equations and thus change to the continuous,  $N\to\infty$, limit. Assuming that all the power lines in the system are of the same grade,  we set $r_k/x_k$ to a constant. Then $r_k = r l_k/L, x_k = x l_k/L$, where $L$ is the total length of the feeder line.
To develop the standard homogenization approach we assume that $P_k, Q_k, v_k$ are all decomposable into a sum of two components, one large but varying smoothly, and another small by amplitude but varying fast:
$F_k = F(z) + \tilde F(L_k)/N$,
where $z=L_k/L$, and $L_k = \sum_{i=0}^{k-1} l_k $; and both $F$ and $\tilde F$ are $O(1)$. The homogenization approach suggests that,  given homogeneous consumption along the feeder, $p_k$ and $q_k$ are small while  $ p(z) = p_k L /l_{k}$ and $q(z) = q_k L/l_k$ are $O(1)$. Relating differences to derivatives,
$ F_{k+1} - F_k \approx  F'(z)l_k/L$,
and then applying it to all the relative differences in the Distflow equations one arrives at the following set of Ordinary Differential Equations (ODEs) that will be called, naturally, DistFlow ODEs
\begin{eqnarray}
\frac{d}{dz}\left(\begin{array}{c}
P\\ Q\\v
\end{array}\right)=
\left(\begin{array}{c}
p-r\frac{P^2+Q^2}{v^2}\\ q-x\frac{P^2+Q^2}{v^2}\\ -\frac{rP+xQ}{v}
\end{array}\right).
\label{DistFlow-ODE}
\end{eqnarray}
The DistFlow ODEs need to be solved under the following mixed
(as defined on both ends of the feeder) conditions
\begin{eqnarray}
v(0)=1,\quad P(L)=Q(L)=0. \label{BC}
\end{eqnarray}

\section{Inverter based control of reactive power}
\label{sec:Control}
Many nodes of existing feeders are capable to control their reactive consumption/production. This freedom is particularly handy in the modern PV systems which contain inverters \cite{lopes2007integrating}. Not only these inverters are capable of controlling the reactive power (obviously within their natural limits,  usually related to the maximum level of real power the inverter is calibrated for) but they are also capable of doing so (at least in principle) in a smart and almost instantaneous manner, e.g. responding to some other measurements available for the feedback. A number of inverter-based control schemes responding to local voltage were recently discussed \cite{lopes2007integrating,EPRI2010}. Here,  we will limit ourselves to analysis of the following two simple control schemes, which should all be contrasted against the ``do nothing" case,  when reactive power is injected by inverter in a natural proportion to the respective real power injection.  Our first model (and very simple) control case is the one of the ``zero power factor",  when the reactive power flow from the node in the feeder is simply set to zero, $q=0$. Another,  and a bit more sophisticated, model analyzed here is similar to the EPRI example from \cite{EPRI2010}, also discussed in \cite{10TSBCa,10TSBCb,11TSBC}:
\begin{eqnarray}
q(v)=q_0\left(1-\frac{2}{1+\exp(-4(v-1)/\delta)}\right),
\label{epri}
\end{eqnarray}
where $q_0$ and $\delta$ are positive constants. $q_0$ is related to control capacity of the inverter and $\delta$ is the parameter expressing tolerance in the level of voltage variations. In this scheme with local feedback one assumes that $q(z)$ responds to voltage measured at the same location, $v(z)$.

\section{Re-scaling and Simplification of the Boundary Value ODE Problem}
\label{sec:Re-scaling}
Assuming that $p=\mbox{const}$, and changing from the flow densities $P,Q$, voltage, $v$, and position along the feeder, $z$, to the new variables
\begin{eqnarray}
&& \rho(s)=\sqrt{\frac{r}{|p|}}\frac{P(z)}{v(L)},\quad
\tau(s)=\sqrt{\frac{r}{|p|}}\frac{Q(z)}{v(L)},\nonumber\\
&& \upsilon(s)=\frac{v(z)}{v(L)},\quad
s=\frac{\sqrt{|p|r}}{v(L)}(L-z),
\label{new-var}
\end{eqnarray}
one arrives at the following rescaled version of the DistFlow-ODEs (\ref{DistFlow-ODE})
\begin{eqnarray}
-\frac{d}{d s}\left(\begin{array}{c}
\rho\\ \tau\\ \upsilon
\end{array}\right)=
\left(\begin{array}{c}
\mbox{sign}(p)-\frac{\rho^2+\tau^2}{\upsilon^2}\\
A-B\frac{\rho^2+\tau^2}{\upsilon^2}\\
-\frac{\rho+B\tau}{\upsilon}
\end{array}\right),
\label{DistFlow-ODE-new}
\end{eqnarray}
where $A\equiv q/|p|, B\equiv x/r$ and the boundary conditions (\ref{BC}) turn into
\begin{eqnarray}
\upsilon(0)=1,\quad \rho(0)=\tau(0)=0.
\label{BC-new}
\end{eqnarray}

The rescaled formulation is advantageous as
when $q=\mbox{const}$ the resulting system of Eqs.~(\ref{DistFlow-ODE-new},\ref{BC-new}) defines a Cauchy problem, with the conditions (\ref{BC-new}) stated only on one end of the $s$, interval, in a pleasing contrast with the mixed original formulation of (\ref{DistFlow-ODE},\ref{BC}). We solve the Cauchy problem (\ref{DistFlow-ODE-new},\ref{BC-new}) forward in the rescaled time, $s: 0\to s_*$,  thus arriving at $\rho(s_*), \tau(s_*), \upsilon(s_*)$, and then one recomputes $L$, as well as $P(0),Q(0)$ and $v(L)$, according to
\begin{eqnarray}
&& L=\frac{s_*}{\upsilon(s_*)\sqrt{|p| r}},\quad v(L)=\frac{1}{\upsilon(s_*)},
\nonumber\\
&& P(0)=\frac{\rho(s_*)\sqrt{|p|/r}}{\upsilon(s_*)},\quad
Q(0)=\frac{\tau(s_*)\sqrt{|p|/r}}{\upsilon(s_*)}.
\label{recomp}
\end{eqnarray}
To summarize, Eqs.~(\ref{DistFlow-ODE-new},\ref{BC-new},\ref{recomp}) allow efficient computations of the original mixed problem (\ref{DistFlow-ODE},\ref{BC}) for different values of the feeder length $L$ by simply scanning (increasing) $s_*$.

Note that the trick of  reducing the original mixed problem  Eqs.~(\ref{DistFlow-ODE-new},\ref{BC-new},\ref{recomp}) to single swap Cauchy computations generalizes to the special cases of the zero power factor control, $q=0$. However,  the technique does not apply to the case of a voltage dependent reactive control. Therefore, in discussing results of our simulations in the next Section we relay on the efficient Cauchy technique,  but also (and when the latter does not apply) on resolving the boundary value problem directly.

\section{Case Studies of the DistFlow ODEs and Discussions}
\label{sec:Simulations}
This Section is devoted to analysis of the DistFlow ODEs in a number of interesting regimes. Utilizing the rescaling arguments above, one sets $r=x=1$ in all of the simulations. We consider only lines constructed of statistically similar customers. Therefore,  and in view of the homogenization and scaling arguments above, we will set the level of real power to $p=\mp 1$, depending on if the feeder is in the consumption or production regime. When no extra voltage control applies, this constant power consumption/production also translates into a constant consumption/production of the reactive power. We will set the base case of constant, $q=p/2$ reactive consumption/production and than compare it with two simple control schemes (assumed realizable by distributed inverters): zero power-factor control, $q=0$, and feedback control on voltage described by Eq.~(\ref{epri}).

\begin{figure}
 \begin{center}
 \includegraphics[width=0.5\textwidth]{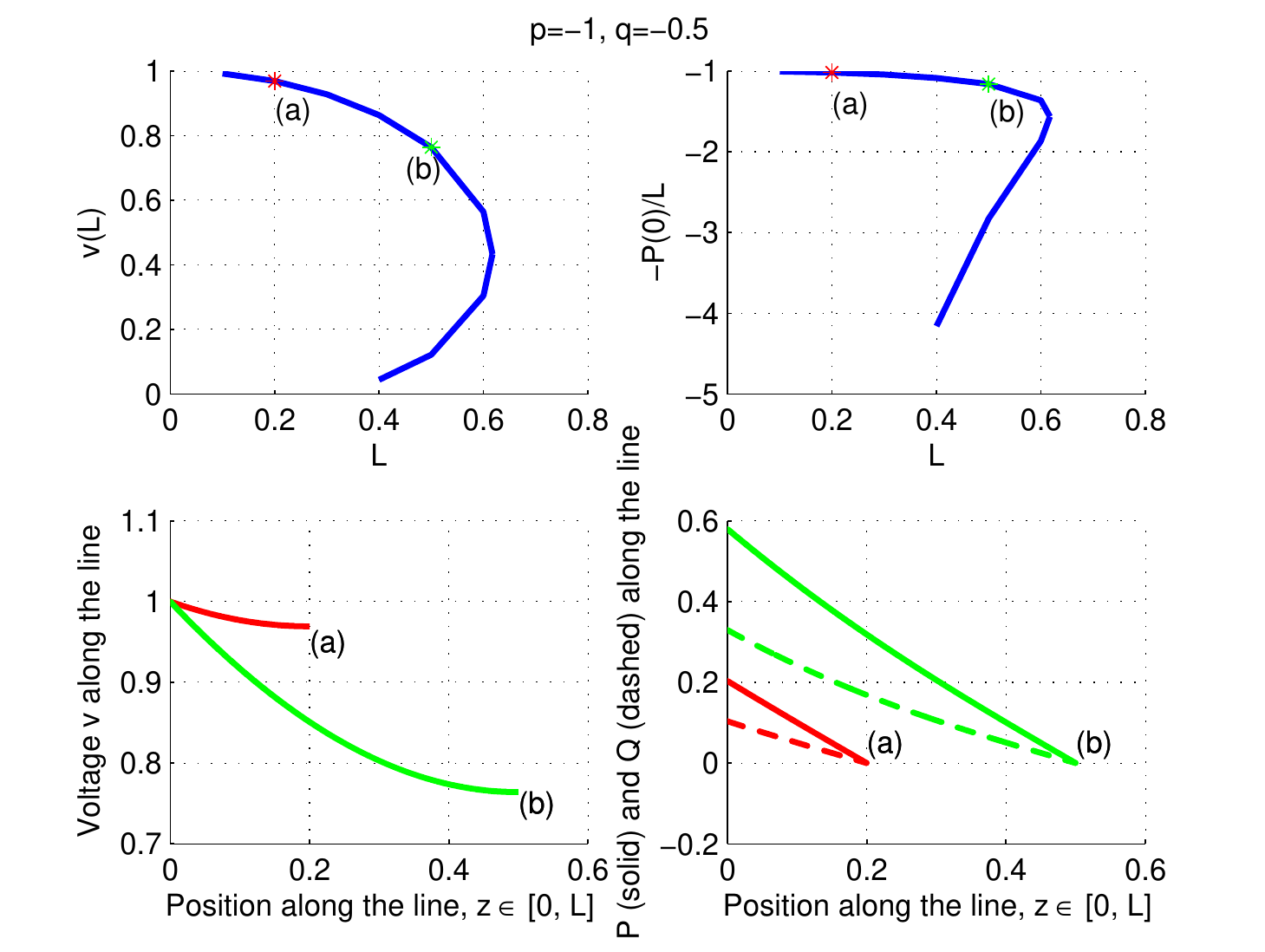}
 \caption{The case of uniform distributed consumption of real and reactive powers.}
 \label{fig:1}
 \end{center}
\end{figure}

To set up the stage we,  first, analyze the standard case,  where both $p$ and $q$ are constant which are negative along the feeder, thus describing the feeder build of customers with constant (voltage independent) consumption. To generate data shown in Fig.~(\ref{fig:1}) we fix $p=-1$ and $q=-0.5$ and then scan different values of $L$. The result, shown  on the top left sub-Figure of Fig.~(\ref{fig:1}), is a standard in power engineering nose-curve illustrating the fact that the amount of power drawn from the system cannot exceed the threshold dependent on the system characteristics.  Figure on the top right of Fig.~(\ref{fig:1}) shows trade-off between the power-utilization of the feeder, measured in the ratio of the power injected at the head of the line, $P(0)$ to the power consumed along the line, $-p*L=L$,
and the voltage drop at the end of the line, $v(L)$.
Two bottom sub-figures in Fig.~(\ref{fig:1}) show profiles of voltage and power flows for lines with parameters correspondent to green and red markers in the top sub-figures. Markers select two representative stable solutions. (Here stability is defined with respect to the standard exciter dynamics, responding to the change in voltage. Then, the upper portion of the nose curve in Fig.~(\ref{fig:1}), where $dv(L)/dP(L)>0$, is stable, while its lower portion,  were $dv(L)/dP(L)<0$, is unstable. We note that this definition of stability is a proxy to actual stability,  which should be analyzed in the future in reference to actual dynamical equations,  which are in fact different for different types of consumers/producers.)

\begin{figure}
 \begin{center}
 \includegraphics[width=0.5\textwidth]{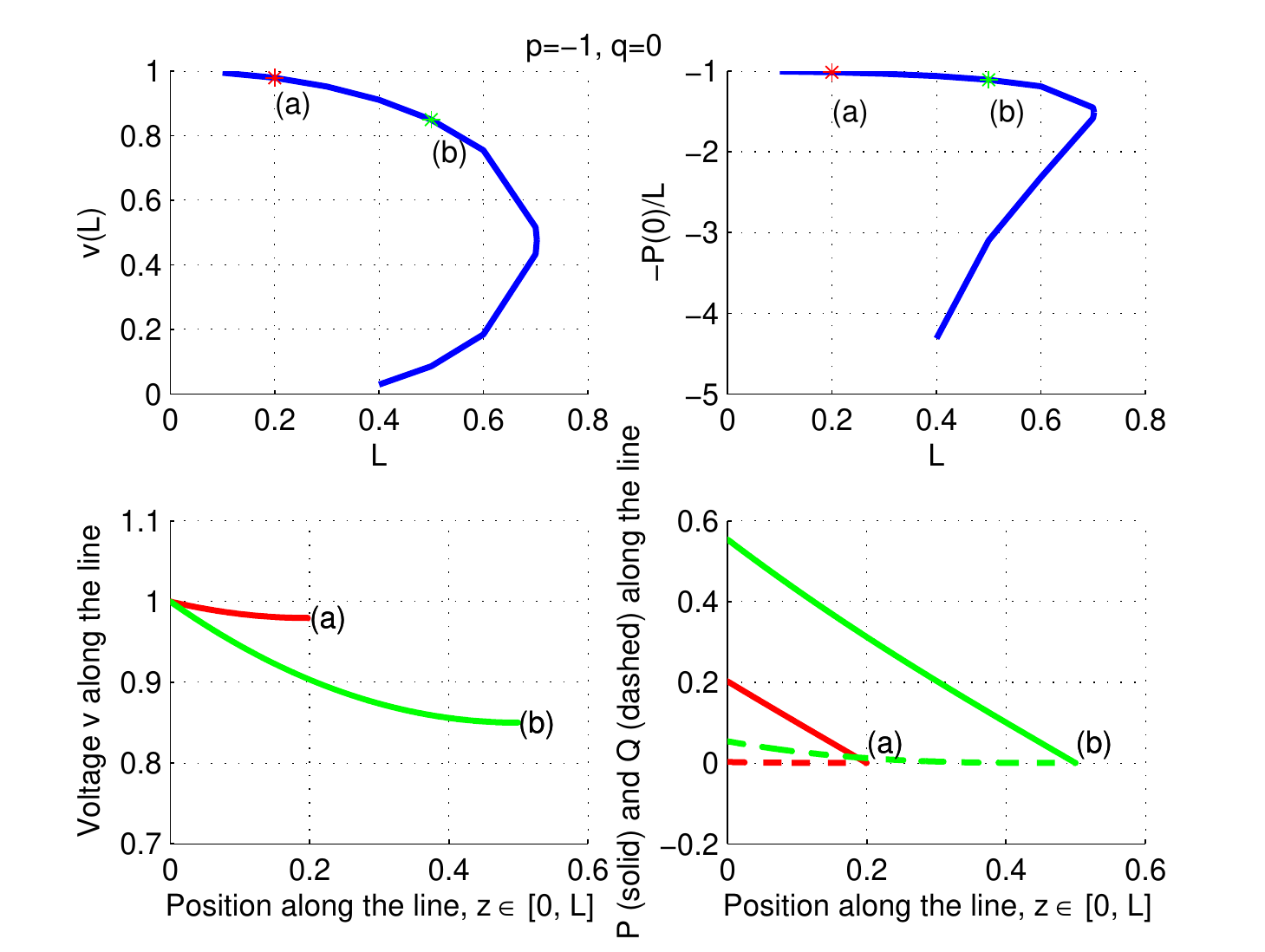}
 \caption{The case of uniform consumption and zero-power factor
 reactive control.}
 \label{fig:2}
 \end{center}
\end{figure}

\begin{figure}
 \begin{center}
 \includegraphics[width=0.5\textwidth]{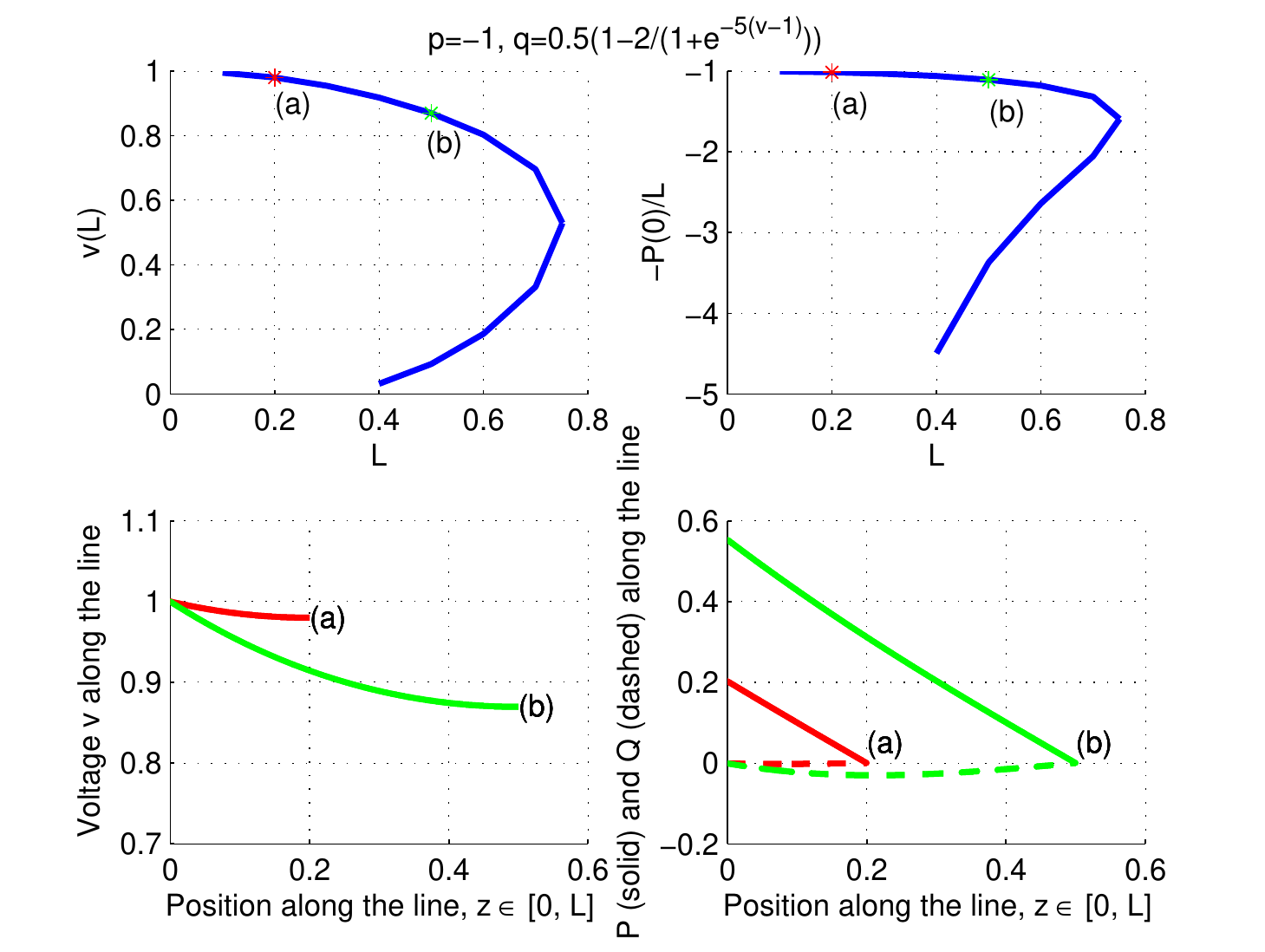}
 \caption{The case of uniform consumption and reactive control with feedback loop on voltage.}
 \label{fig:3}
 \end{center}
\end{figure}

The bare case of distributed consumption,  shown in Fig.~(\ref{fig:1}), should be compared with Fig.~(\ref{fig:2}) and
Fig.~(\ref{fig:3}). We use the same type of diagnostics for creating the complex Figures as before, however we apply it now, in Fig.~(\ref{fig:2}) and Fig.~(\ref{fig:3}), to the case with the zero power-factor control ($q=0$) and the reactive control of Eq.~(\ref{epri}) with feedback on voltage. Changing from the bare case to the zero  power factor control does improve characteristics of the system by shifting the nose curve to the right and thus reducing voltage decrease for the feeder of the fixed length. The improvement continues with further transition from the zero power factor case of  Fig.~(\ref{fig:2}) to the case of reactive feedback control of Fig.~(\ref{fig:3}). These effects are quantitative and intuitive,  and in fact not surprising as also following from an even simpler model with the  distributed feeder replaced by a single clamped/aggregated load.

\begin{figure}
 \begin{center}
 \includegraphics[width=0.5\textwidth]{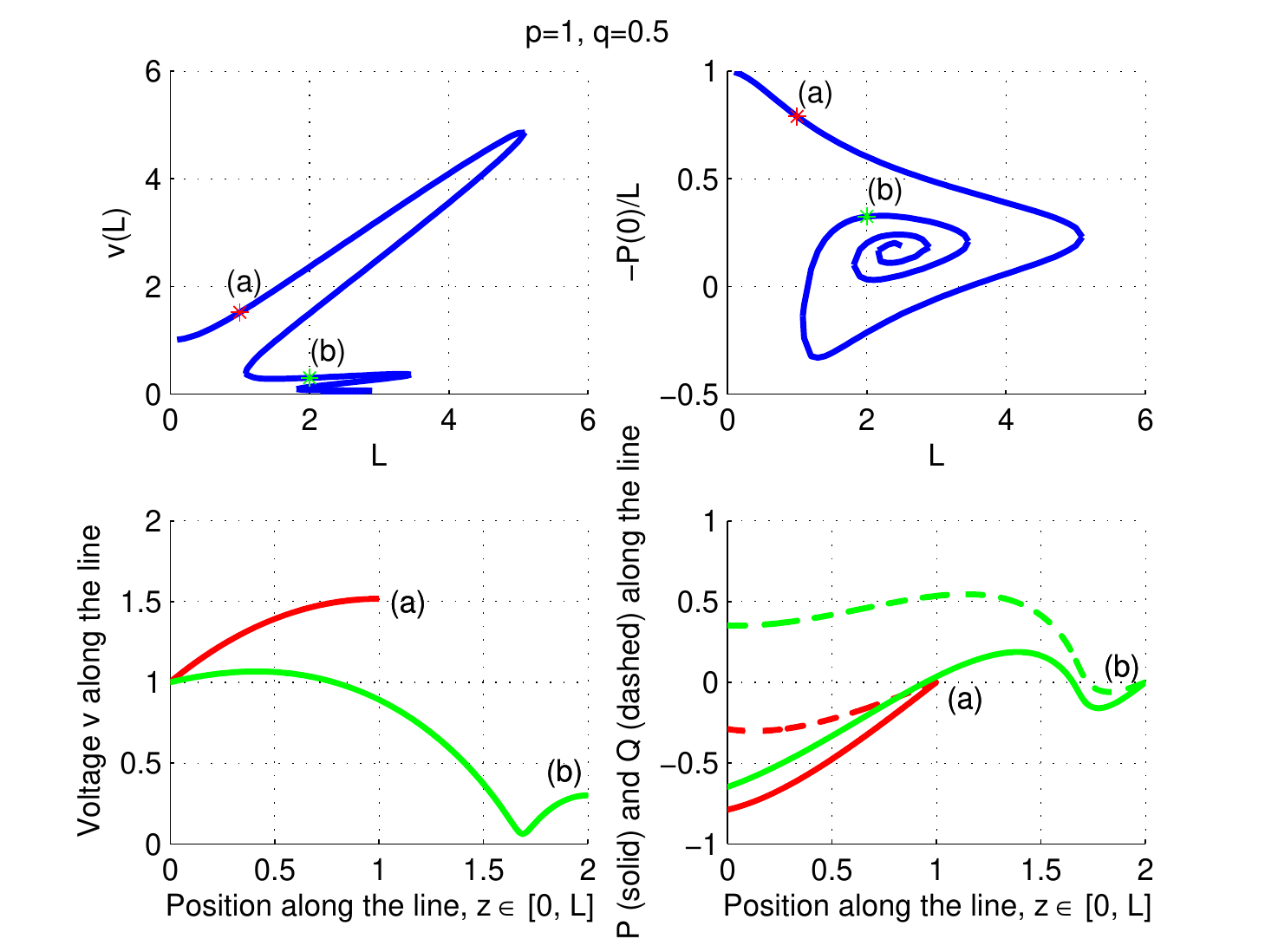}
 \caption{The case of uniformly distributed and comparable generation of real and reactive powers}
 \label{fig:4}
 \end{center}
\end{figure}

\begin{figure}
 \begin{center}
 \includegraphics[width=0.5\textwidth]{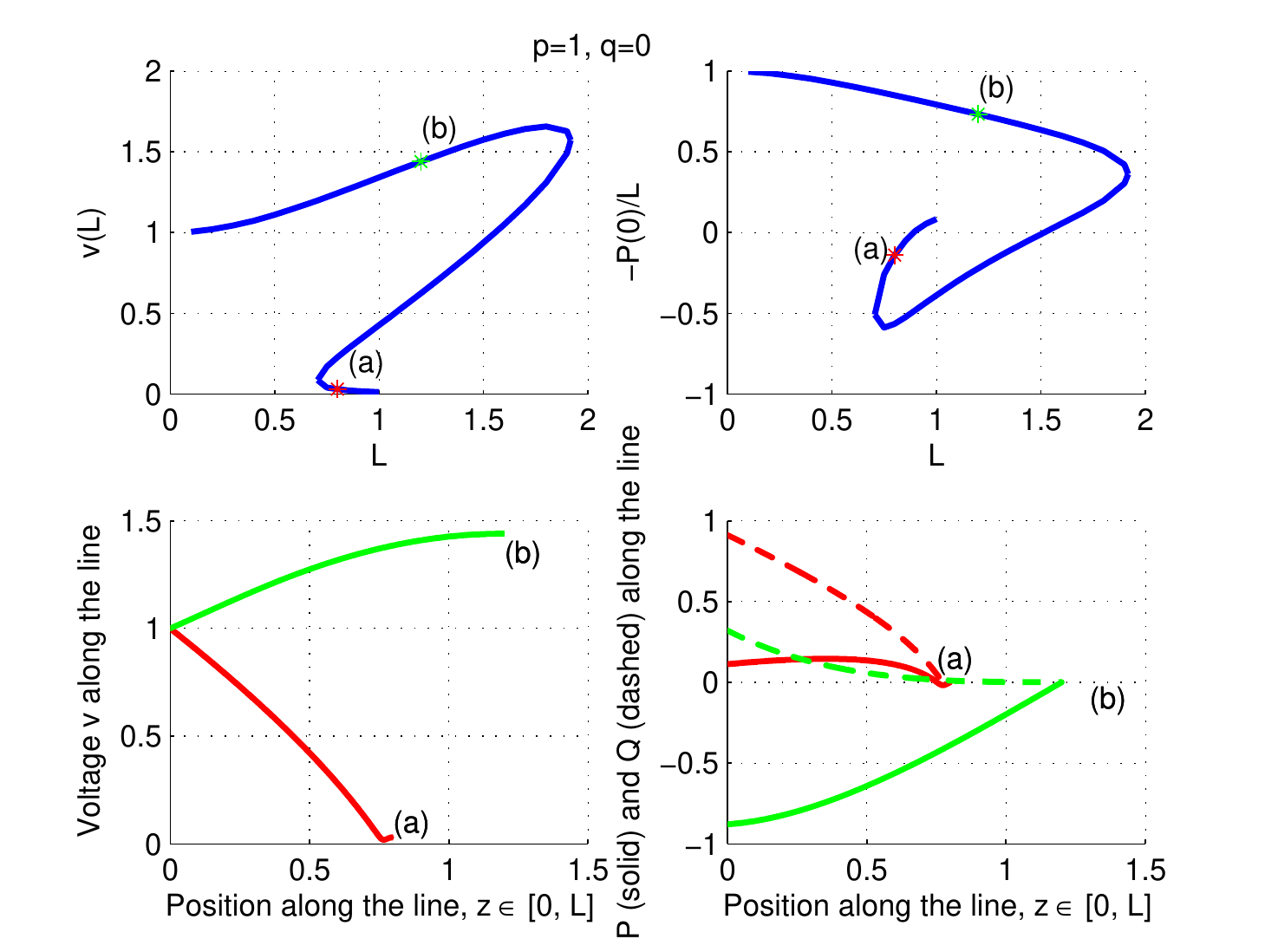}
 \caption{The case of uniformly distributed generation of real power and zero reactive power (zero power factor) control}
 \label{fig:5}
 \end{center}
\end{figure}

\begin{figure}
 \begin{center}
 \includegraphics[width=0.5\textwidth]{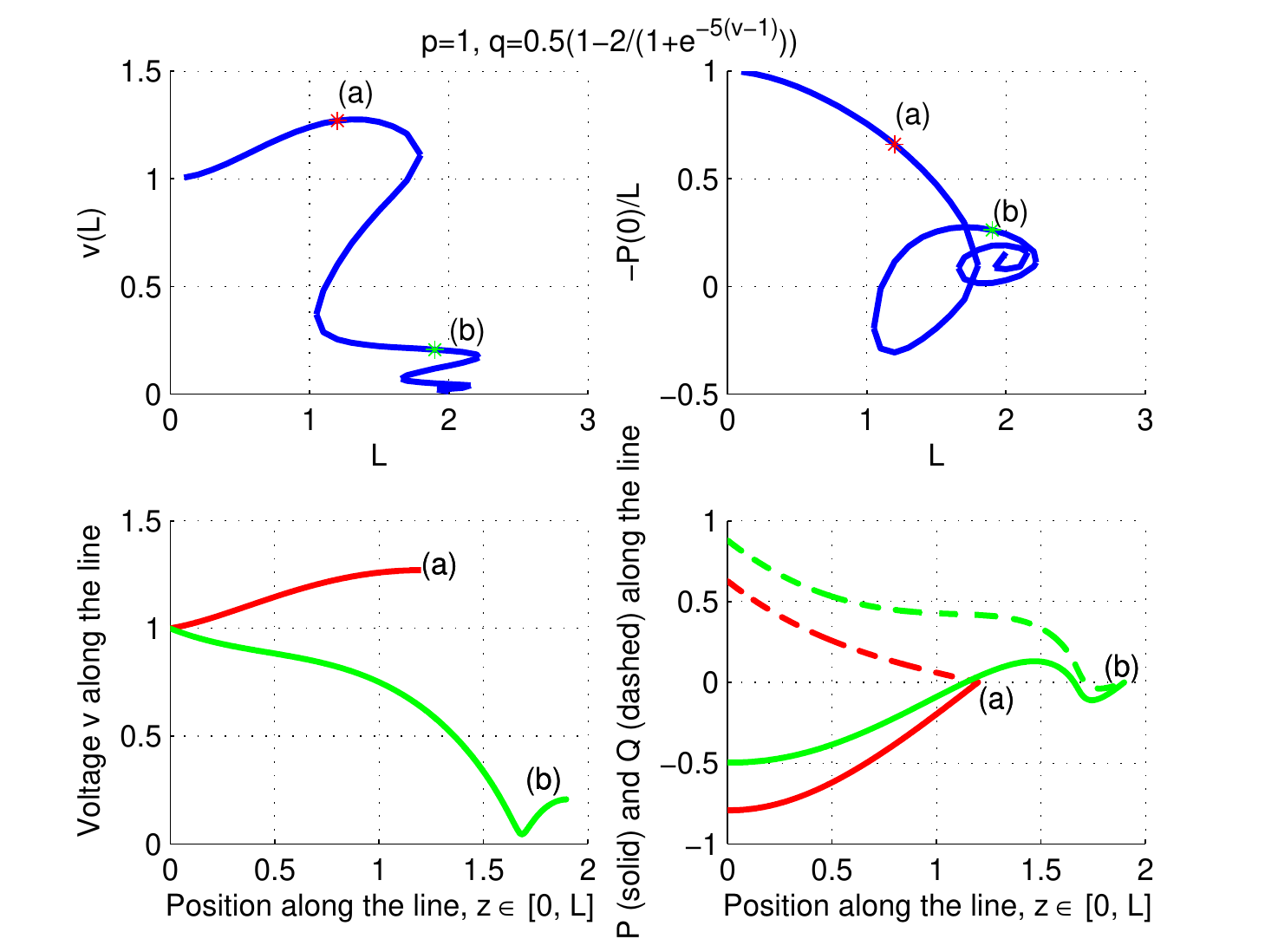}
 \caption{The case of uniformly distributed generation of real power and reactive control with feedback on voltage.}
 \label{fig:6}
 \end{center}
\end{figure}

Next we turn to analysis of a feeder with distributed generation, which may be though of as a homogenized proxy for a feeder rich on  Photo-Voltaic systems generating enough of power during a sunny day to compensate for local consumption and injecting the access of power into the system. As before we showed in Fig.~(\ref{fig:4},\ref{fig:5},\ref{fig:6}) the bare case (of comparable injection of real and reactive powers), the case of the zero power-factor ($q=0$) control,  and the case of reactive control with feedback to voltage respectively.  We use the same diagnostics as above, showing voltage profile and power utilization as functions of the varying length of the feeder in the top sub-figures,  and voltage and power flow profiles analyzed as functions of the position along the feeder, for two (marked by green and red markers) representative cases.

We observe that behavior in the Figures with distributed generation is qualitatively different from the one seen above in the case of distributed consumption. Most interesting new phenomenon,  the one which was not seen in the case of distributed consumption, and which would not be seen in the naive model replacing the whole generating feeder by a single effective node, is associated with the emergence of multiple stable solutions.  Red and green markers on the top sub-Figures of Fig.~(\ref{fig:4}) show solutions
from the top (in terms of voltage) and the second top stable branches. Even though the lower voltage solution may seem unreachable as showing dangerously low voltage, the bare co-existence (for a range of L) of multiple stable solutions itself is extremely dangerous for the distribution system operations.  Indeed, fault in the transmission system, resulting in an abrupt degradation  of voltage to zero at the head of the line, may force the system to equilibrate in the post-fault regime to an undesirable low voltage branch. This scenario would be similar to the one associated with the phenomenon of the Fault Induced Delayed Voltage Recovery (FIDVR) \cite{Hill1993,Davy1997} which was observed in feeders rich on inductive motors (e.g. many air-condition loads with small inertia) \cite{Lesieutre2005}. The phenomenon of FIDVR can be interpreted as the one coming from a complicated hysteretic $p(v)$ and $q(v)$ characteristics of the individual inductive loads \cite{12BC}.   Similarly,  we should worry about a hysteretic behavior in the case of a feeder with distributed generation.  Notice, however that here the phenomenon  does not originate from  individual loads,  but it rather emerges as a result of a complex spatial distribution of voltage and power flows along the feeder.

We analyze details of voltage and power flow profiles for two representative feeders, selected in the top
portion of Fig.~(\ref{fig:4}) by green and red markers. For the top branch (marked green) the profiles are standard:  voltage raises monotonically,  both real and reactive power flows are injected on the top of the line and then decrease monotonically (by absolute values) as we move along the line. However, the profiles are distinctly different in the second regime (marked by red markers in the top sub-Figures). Here, the direction of real flow reverses in the middle of the feeder - instead of flowing towards the head of the feeder (which would be the naively anticipated case with the export of power from this part of the feeder) the flow changes to local import,  with the power flowing towards the end of the feeder. Note that this oscillatory behavior may be even more complex (with multiple reversals of the flow direction along the feeder) for other lower voltage branches.

The two cases of control,  applied to this bare case with distributed generators are shown in Fig.~(\ref{fig:5}) and Fig.~(\ref{fig:6}).  The changes in these cases (when compared with the bare case) are visible,  however not alternating the principle observation made above: multiple stable solutions may still be realized in the feeder. We can understand this observation as follows.  The main design suggestion for these schemes were based on considerations of voltage control.  In this regards,  the controls do achieve the anticipated goal as they do reduce increase of voltage in the top,  base solution.  However,  the controls do not change the potentially dangerous co-existence of the reasonable (in terms of its voltage and power flow profiles) solution with its complex and very undesirable low-voltage sibling(s).

One may question validity of modeling the distributed generation by generation with output independent of voltage in the low-voltage regime.  We have tested a more realistic case (regularized at low voltage)  as well (not shown in Figures) and observe that even if $p$ decreases to zero with $v$ in some reasonable manner,  the main conclusion drawn above within the constant generation model stays - multiple stable solutions may emerge and standard reactive control, taking care of the voltage problems in the top branch,  does not remove the operationally dangerous multiplicity/degeneracy of the power flow solutions.

\section{Conclusions and Path Forward}
\label{sec:Conclusions}
The key findings of this paper can be summarized as follows:
\begin{itemize}
	\item We have presented a novel approach to analysis of a long feeder with statistically homogeneous distribution of loads. When the number of consumers/buses is large, equations describing power flows over the feeder can be reduced to a system of ordinary differential equations. The approach is computationally advantageous as containing fewer number of parameters than in the original discrete formulation. The technique also allows to account for distributed reactive and active power controls of individual loads.
	\item In the (smart grid) future we expect some feeders to perform in the new regime of distributed generation. We show that this regime exhibits a number of features distinctly different from the standard feeder of today containing only distributed consumption. We observe that in the distributed generation regime the power flow may reverse its direction along the feeder multiple number of times. Moreover, under certain conditions bi- and multi-stability can be observed, thus resulting in a hard to control,  and dangerous, hysteretic behavior.
	\item The distributed control of reactive power that is possible with programmable inverters provides principal capability of modifying the solution manifold. This option can be used to recover from an abrupt voltage fault as well as for trading off between losses and voltage deviations in the system.
\end{itemize}

Note that the model(s) considered in the manuscript represent only rough qualitative proxies for real power systems with distributed generation. In spite of these models simplicity they show rich and nontrivial solutions,  thus providing a short glimpse into the emerging complexity associated with the large-scale integration of distributed generation in modern power systems. More realistic and detailed models of power distribution systems should account for non-trivial voltage dependence of aggregated loads and generators on every bus. Accurate modeling is especially important for describing recovery from the transmission level induced faults in the feeder which was operating in the regime of distributed generation prior to the fault. It will also be important in the future to have a detailed dynamic analysis of the post-fault transients, in particular accounting for intermittent fluctuations of the distributed renewable sources. Finally, we anticipate that the technique introduced in this manuscript will be instrumental for addressing the fundamentally important (but still widely open) questions of formulating and designing optimal, centralized or distributed, control of reactive and active power flows in the power distribution systems.

\section{Acknowledgements}

We are thankful to the participants of the ``Optimization and Control for Smart Grids" LDRD DR project at Los Alamos and Smart Grid Seminar Series at CNLS/LANL, and especially to S. Backhaus for multiple fruitful discussions. 

\bibliographystyle{IEEETran}
\bibliography{voltage}

\end{document}